\documentclass[12pt]{elsarticle}
\usepackage{amssymb,amsmath,amsthm,enumerate}
\usepackage{amssymb}

 \theoremstyle{definition}
 
 \theoremstyle{remark}
 
 \theoremstyle{example}
 
\begin{document}
\begin{frontmatter}
\title{A note on Sonnenschein summability matrices}
\author{Gholamreza Talebi*, Masoud Aminizadeh} 
\address{Department of Mathematics,\\ Vali-e-Asr University of Rafsanjan,\\ Rafsanjan, Islamic Republic of Iran.\\E-mail:  Gh.talebi@vru.ac.ir,\vspace{-1cm}}
\begin{abstract}
In this note, we give a simple method for computing the column sums of the Sonnenschein summability matrices. 
\end{abstract}
\end{frontmatter}
\section{Introduction}
Let $f(z)$ be an analytic function in  $|z|<r, r\geq1$ with $f(1)=1.$ The matrix $S=(f_{n,k}),$ where $(f_{n,k})$  are defined by ${\left[ {f\left( z \right)} \right]^n} = \sum_{k = 0}^\infty  {{f_{n,k}}{z^k}}$ is called a Sonnenschein matrix \cite{J}. The special choice 
\[ f(z)=\frac{\alpha+\left(1-\alpha-\beta\right)z}{1-\beta z},~z\in\mathbb{C}\backslash\{\frac{1}{\beta}\},\]
where $\alpha$ and $\beta$ are complex numbers, gives the Karamata matrix $K[\alpha,\beta]$ and its coefficients as a Sonnenschein matrix are given by \cite{Boos}
\[{f_{n,k}} = \sum\limits_{v = 0}^k {\left( \begin{array}{l}
n\\
v
\end{array} \right){{\left( {1 - \alpha  - \beta } \right)}^v}{\alpha ^{n - v}}\left( \begin{array}{l}
n + k - v - 1\\
\,\,\,\,\,\,\,\,\,\,k - v
\end{array} \right)
{\beta ^{k - v}}}.\]\par
Recently in \cite{WALA}, the authors have calculated the row and column sums of Karamata matrices in a relatively complicated way. In this note we give a new and simple method for computing the column sums of these matrices which can be also applied to another Sonnenschein matrices. \vspace{.2cm}\par
Start with a general $f.$ It is clear that if $|f(0)|<1$ then $|f(z)|<1$ in a neighborhood of $z=0.$ Thus
\[\sum\limits_{n = 0}^\infty  {{{\left[ {f\left( z \right)} \right]}^n}}  = \frac{1}{{1 - f(z)}}.\]
Now the sum we want is the coefficient of $z^k$ of the right-hand side of the above equation. For the case of the Karamata matrices we have $f(0)=\alpha,$ so we assume first that $|\alpha|<1.$ Now 
\[\frac{1}{{1 - f(z)}} = \frac{1}{{1 - \alpha }} + \frac{{\left( {1 - \beta } \right)z}}{{\left( {1 - \alpha } \right)\left( {1 - z} \right)}}=\frac{1}{{1 - \alpha }} + \frac{{1 - \beta }}{{1 - \alpha }}\sum\limits_{k = 1}^\infty  {{z^k}} .\]
The coefficient of $z^k$ is then easily found. Indeed, the sum of the first column is $\frac{1}{{1 - \alpha }} $, and the sum of all other columns are $ \frac{{1 - \beta }}{{1 - \alpha }}.$ \vspace{.2cm}\par
The point is that this method may apply to other choices of the function $f$ and other Sonnenschein matrices. 
As another example consider the function 
\[h(z) = sin^2(\frac{\pi z}{2}),~~(z\in\mathbb{C})\]
which is holomorphic function and its coefficients as a Sonnenschein matrix are as followings:
\begin{enumerate}
\item If $n=0$, then ${a_{0,k}} =\delta_{0k},$ for all $k;$
\item If $n\neq 0$ and $k=0$, then ${a_{n,k}} =\frac{1}{{{4^n}}}\left( {_{\,\,n}^{2n}} \right) + \sum\limits_{r = 0}^{n - 1} {\frac{{{{( - 1)}^{n + r}}}}{{{2^{2n - 1}}}}\left( {_{2r}^{2n}} \right)}$;
\item  If $n\neq 0 $ and $k\neq~0$, then 
$~~{a_{n,k}} = \left\{ \begin{array}{l}
\,\,\,\,0 \,\,\,\,\,\,\,\,\,\,\,\,\,\,\,\,\,\,\,\,\,\,\,\,\,\,\,\,\,\,\,\,\,\,\,\,\,\,\,\,\,\,\,\,\,\,\,\,\,\,\,\,\,\,\,\,\,\,\,\,\,\,\,\,\,\,\, k$ is odd$,\\
\\
{\pi ^{2k}}\sum\limits_{r = 0}^{n - 1} {\frac{{{{\left( { - 1} \right)}^{n + r + k}}{{\left( {n - r} \right)}^{2k}}}}{{{2^{2n - 1}}\left( {2k} \right)!}}\left( {_{2r}^{2n}} \right)} \,\,\,\,\,\,\,\,k $ is even.$
\end{array} \right.$
\end{enumerate}
We have $h(0)=0$, thus 
\[\begin{array}{l}
\frac{1}{{1 - h(z)}} = {\sec ^2}\left( {\frac{{\pi z}}{2}} \right) =\frac{2}{\pi} \frac{{d\tan \left( {\frac{{\pi z}}{2}} \right)}}{{dz}}\\
\\
\,\,\,\,\,\,\,\,\,\,\,\,\,\,\,\,\,\,\,\,\,\,\,\,\,\,\,\,\,\,\,\,\,\,\,\,\,\,\,\,\,\,\,\,\,\,\,\,\,=\frac{2}{\pi} \frac{{d\left( {\sum\limits_{n = 0}^\infty  {\frac{{{{( - 1)}^n}2({2^{2n + 2}} - 1)}}{{\left( {2n + 2} \right)!}}{B_{2n + 2}}~{\pi ^{2n + 1}}} {z^{2n + 1}}} \right)}}{{dz}}\,\\
\\
\,\,\,\,\,\,\,\,\,\,\,\,\,\,\,\,\,\,\,\,\,\,\,\,\,\,\,\,\,\,\,\,\,\,\,\,\,\,\,\,\,\,\,\,\,\,\,\,\,= \sum\limits_{n = 0}^\infty  {\frac{{{{( - 1)}^n}\left( {8n + 4} \right)({2^{2n + 2}} - 1)}}{{\left( {2n + 2} \right)!}}{B_{2n + 2}}~{\pi ^{2n }}} {z^{2n}},
\end{array}\]
where $B_n$ is the sequence of Bernoulli numbers (see \cite{Po}, pp. 274 - 275), defined by 
\[{B_n} = \sum\limits_{k = 0}^n {\frac{1}{{k + 1}}\sum\limits_{r = 0}^k {{{( - 1)}^r}\left( {_r^k} \right)} } {r^n}.~~~~ (n=0,1,2,...)\]
Therefore, the sum of odd columns are $0$, and the sum of $2n-$th columns ($n=0,1,2,...$) is \[{\frac{{{{( - 1)}^n}\left( {8n + 4} \right)({2^{2n + 2}} - 1)}}{{\left( {2n + 2} \right)!}}{B_{2n + 2}}~{\pi ^{2n}}}.\]
\section*{Acknowledgements}
The authors would like to thank professor Mourad Ismail [University of Central Florida, Orlando] for his technical assistance during the preparation of the this note.\\\\
\textbf{References}
\bibliographystyle{amsplain}

\end{document}